\documentclass[a4paper,twoside,11pt]{article}
\usepackage{amsfonts}
\usepackage{amssymb}\usepackage{amsmath}
\newtheorem{thm}{Theorem}[section]
 
 \newtheorem{lem}[thm]{Lemma}

 \numberwithin{equation}{section}
\newcommand{\double}{\baselineskip 1.24 \baselineskip}

\title{Some new Inequalities of Hardy-Hilbert Type with general kernel
}

\author{{Guang-Sheng  Chen\thanks{\text{E-mail address}: cgswavelets@126.com(Chen)
}\quad}\\
{\small Department of Computer Engineering, Guangxi Modern
Vocational Technology College,} \\{\small Hechi,Guangxi, 547000,
P.R. China}
}
\begin{document}
\date{}
\maketitle \double

\textbf{Abstract:}\quad In this paper, by using hardy inequality, we establish some new
integral inequalities of Hardy-Hilbert type with general kernel. As
applications, equivalent forms and some particular results are built; the
corresponding to the double series inequalities are given.reverse forms are
considered also. 

\section{Introduction}
\hskip\parindent

 If $p > 1$,$1 / p + 1 / q = 1$, $f(x)$, $g(x) \ge 0$, $0 < \int_0^\infty
{f^p(x)dx} < \infty $, and

$0 < \int_0^\infty {g^q(x)dx} < \infty ,$then the famous Hardy--Hilbert's
inequality (see [1]) and an equivalent form are given by
\begin{equation*}
\int_0^\infty {\int_0^\infty {\frac{f(x)g(y)}{x + y}dxdy} } < \frac{\pi
}{\sin (\pi / p)}[\int_0^\infty {f^p(x)dx} ]^{1 / p}[\int_0^\infty
{g^q(x)dx} ]^{\frac{1}{q}},
\tag{1.1}
\label{1.1}
\end{equation*}
And
\begin{equation*}
\int_0^\infty {[\int_0^\infty {\frac{f(x)}{x + y}dx]^pdy} } < [\frac{\pi
}{\sin (\pi / p)}]^p\int_0^\infty {f^p(x)} dx,
\tag{1.2}
\label{1.2}
\end{equation*}
Where the constant factor $\frac{\pi }{\sin (\pi / p)}$ and $[\frac{\pi
}{\sin (\pi / p)}]^p$ are the best possible.

Hardy et al. [1] gave an inequality and its equivalent form , under the same
condition of (1.1), similar to (1.1) as :
\begin{equation*}
\int_0^\infty {\int_0^\infty {\frac{f(x)g(y)}{\max \{x,y\}}dxdy} } <
pq[\int_0^\infty {f^p(x)dx} ]^{1 / p}[\int_0^\infty {g^q(x)dx}
]^{\frac{1}{q}},
\tag{1.3}
\label{1.3}
\end{equation*}
And
\begin{equation*}
\int_0^\infty {[\int_0^\infty {\frac{f(x)}{\max \{x,y\}}dx]^pdy} } <
[pq]^p\int_0^\infty {f^p(x)} dx,
\tag{1.4}
\label{1.4}
\end{equation*}
Where the constant factor $pq$ and $[pq]^p$ are the best possible.

Inequalities (1.1), (1.2) ,(1.3) and (1.4) are important in analysis and its
applications(see [2]). In the recent years , many generalization and refinements of
these inequalities have been also obtained ( see [3-8] ).

Recently Das and Sahoo [8] have given a new inequality similar to
Hardy--Hilbert inequality (1.1) as follows:

Let $p > 1$, $1 / p + 1 / q = 1$, $r,s,\lambda > 0$, $r + s = \lambda
$，$f(x)$, $g(x) \ge 0$,

$F(x) = \int_0^x {f(t)dt} ,
G(x) = \int_0^x {g(t)dt} , $ if $0 < \int_0^\infty {f^p(x)dx} < \infty $,

$0 < \int_0^\infty {g^q(x)dx} < \infty ,$then the following two integral
inequalities holds:
\begin{equation*}
\int_0^\infty {\int_0^\infty {\frac{x^{r - \textstyle{1 \over q} - 1}y^{s -
\textstyle{1 \over p} - 1}F(x)G(y)}{(x + y)^\lambda }dxdy} } <
pqB(r,s)[\int_0^\infty {f^p(x)dx} ]^{1 / p}[\int_0^\infty {g^q(x)dx}
]^{\frac{1}{q}},
\tag{1.5}
\label{1.5}
\end{equation*}
And
\begin{equation*}
\int_0^\infty {[\int_0^\infty {\frac{x^{r - \textstyle{1 \over q} - 1}y^{s -
\textstyle{1 \over p}}F(x)}{(x + y)^\lambda }dx]^pdy} } <
[qB(r,s)]^p\int_0^\infty {f^p(x)} dx,
\tag{1.6}
\label{1.6}
\end{equation*}
where the constant factors $pqB(r,s)$ and $[qB(r,s)]^p$ are the best
possible.

Sulaiman [7, Theorem 1] derived a new integral inequality similar to (\ref{eq3})
as follows:

Let $p > 1$,$1 / p + 1 / q = 1$, $p = \lambda - \alpha - 1 > 1$, $p =
\lambda - \beta - 1 > 1$, $\alpha ,\beta > - 1$,$f(x)$, $g(x) \ge 0$, $F(x)
= \int_0^x {f(t)dt} $,$G(x) = \int_0^x {g(t)dt} $, if $0 < \int_0^\infty
{f^p(x)dx} < \infty $，$0 < \int_0^\infty {g^q(x)dx} < \infty $,then the
following two integral inequalities holds:
\begin{equation*}
\begin{split}
&
 \int_0^\infty {\int_0^\infty {\frac{x^{\textstyle{\beta \over
q}}y^{\textstyle{\alpha \over p}}F(x)G(y)}{\max \{x^\lambda ,y^\lambda
\}}dxdy} } \\
 &< \frac{p^{1 - \textstyle{1 \over p}}q^{1 - \textstyle{1 \over q}}}{(\alpha
+ 1)^{\textstyle{1 \over p}}(\beta + 1)^{\textstyle{1 \over q}}(p - 1)(q -
1)}[\int_0^\infty {f^p(x)dx} ]^{1 / p}[\int_0^\infty {g^q(x)dx}
]^{\frac{1}{q}},\\
 \end{split}
 \tag{1.7}
 \label{1.7}
\end{equation*}

in [7],Sulaiman does not prove whether the constant factor is best possible
or not. very recently, Das and Sahoo [9] have given a new generalization of
(1.7). the constant factor is the best possible to prove.

In this paper, we obtain a generalization of the inequalities (1.5) and
(1.7) with general kernel., the constant factor obtained is the best
possible. First we prove the integral version of the inequality and some
particular results. Then we give the discrete analogue of the inequality.
equivalent forms and reverse forms are considered.

\section{Some Lemmas }
\hskip\parindent
We need the following some inequalities, which are well-known as Hardy's
inequalities (cf. Hardy et al. [1]).
\begin{lem}\label{Lemma 2.1.} If $p > 1$,$f(x) \ge 0$ , $F(x) = \int_0^x {f(t)dt} $, and $0 <\int_0^\infty {f^p(x)dx} < \infty $, then
\begin{equation*}
\int_0^\infty {(\frac{F(x)}{x})^pdx} < (\frac{p}{p - 1})^p\int_0^\infty
{f^p(x)dx},
\tag{2.1}
 \label{2.1}
\end{equation*}
unless $f(x) \equiv 0$,The constant is the best possible.
\end{lem}

\begin{lem}\label{Lemma 2.2.} If $0 < p < 1$,$f(x) \ge 0$ , $F(x) = \int_x^\infty {f(t)dt} $,
and $0 < \int_0^\infty {f^p(x)dx} < \infty $, then
\begin{equation*}
\int_0^\infty {(\frac{F(x)}{x})^pdx} > (\frac{p}{1 - p})^p\int_0^\infty
{f^p(x)dx},
\tag{2.2}
 \label{2.2}
\end{equation*}
unless $f(x) \equiv 0$,The constant is the best possible.
\end{lem}

\begin{lem}\label{Lemma 2.3.} If $p > 1$,$a_n \ge 0$, and $A_n = \sum\limits_{i = 1}^n {a{ }_i}
$, then
\begin{equation*}
\sum\limits_{n = 1}^\infty {(\frac{A_n }{n})^p} < (\frac{p}{p -
1})^p\sum\limits_{n = 1}^\infty {a_n^p },
\tag{2.3}
 \label{2.3}
\end{equation*}
unless all the $a_n = 0$. The constant is the best possible.
\end{lem}
\begin{lem}\label{Lemma 2.4.} If $0 < p < 1$,$a_n \ge 0$, and $A_n = \sum\limits_{i = 1}^n {a{
}_i} $, then
\begin{equation*}
\sum\limits_{n = 1}^\infty {(\frac{A_n }{n})^p} < (\frac{p}{1 -
p})^p\sum\limits_{n = 1}^\infty {a_n^p },
\tag{2.4}
 \label{2.4}
\end{equation*}
unless all the $a_n = 0$. The constant is the best possible
\end{lem}
If$k_\lambda (x,y)$ is a measurable function, satisfying for $\lambda $,
$u$, $x$, $y > 0$, $k_\lambda (ux,uy) = u^{ - \lambda }k_\lambda (x,y)$,
then we call $k_\lambda (x,y)$ the homogeneous function of -$\lambda
$-degree.
\begin{lem}\label{Lemma 2.5.} If $r,s,\lambda > 0$,$r + s = \lambda $,$k_\lambda (x,y)
> 0$ is a homogeneous function of -$\lambda $-degree, and $k_\lambda
(r): = \int\limits_0^\infty {k(u,1)u^{r - 1}du} $ a positive number, define
the weight functions $\omega _\lambda (s,x)$ and $\omega _\lambda (r,y)$ as
\begin{equation*}
\omega _\lambda (s,x) = \int_0^\infty {k_\lambda (x,y)x^ry^{s - 1}dy},
\tag{2.5}
\label{2.5}
\end{equation*}
\begin{equation*}
\omega _\lambda (r,y) = \int_0^\infty {k_\lambda (x,y)x^{r - 1}y^sdx},
\tag{2.6}
\label{2.6}
\end{equation*}
then we have

(i) $\int\limits_0^\infty {k(1,u)u^{s - 1}du} = k_\lambda (r)$;

(ii) $\omega _\lambda (s,x) = \omega _\lambda (r,y) = k_\lambda (r)$.
\end{lem}
{\bf Proof.} (i) Setting $v = \frac{1}{u}$, by the assumption, we obtain
\[
\int\limits_0^\infty {k(1,u)u^{s - 1}du} = \int\limits_0^\infty {k(v,1)v^{r
- 1}dv} = k_\lambda (r).
\]

(ii) Setting $u = y / x$ in the integrals $\omega _\lambda (s,x)$, in view
of (i), we still find that$\omega _\lambda (s,x) = k_\lambda (r)$. Similarly
we have $\omega _\lambda (r,y) = k_\lambda (r)$, The lemma is proved.

\begin{lem}\label{Lemma 2.6.} If $p > 1$,$1 / p + 1 / q = 1$, $r,s,\lambda > 0$, $r + s = \lambda
$,$k_\lambda (x,y) > 0$ is a homogeneous function of -$\lambda $-degree,
and $k_\lambda (r): = \int\limits_0^\infty {k(u,1)u^{r - 1}du} $ a positive
number, for sufficiently small $\varepsilon > 0$,setting

\begin{equation*}
I_1 = \int_1^\infty {\int_1^\infty {k_\lambda (x,y)x^{r -
\textstyle{\varepsilon \over p} - 1}y^{s - \textstyle{\varepsilon \over q} -
1}dxdy} } ,
\tag{2.7}
\label{2.7}
\end{equation*}
\begin{equation*}
I_2 = \int_1^\infty {\int_1^\infty {k_\lambda (x,y)x^{r - \textstyle{1 \over
q} - 1}y^{s - \textstyle{\varepsilon \over q} - 1}dxdy} } ),
\tag{2.8}
\label{2.8}
\end{equation*}
\begin{equation*}
I_3 = \int_1^\infty {\int_1^\infty {k_\lambda (x,y)x^{r -
\textstyle{\varepsilon \over p} - 1}y^{s - \textstyle{1 \over p} - 1}dxdy} }
),
\tag{2.9}
\label{2.9}
\end{equation*}
\begin{equation*}
\int_1^\infty {x^{ - 1 - \varepsilon }[\int_0^{1 / x} {k_\lambda (1,u)u^{s -
\textstyle{\varepsilon \over q} - 1}du} ]dx} = O_1 (1),
\tag{2.10}
\label{2.10}
\end{equation*}
\begin{equation*}
\int_1^\infty {x^{ - \textstyle{{\varepsilon + 1} \over q} - 1}[\int_0^{1 /
x} {k_\lambda (1,u)u^{s - \textstyle{\varepsilon \over q} - 1}du} ]dx} = O_2
(1),
\tag{2.11}
\label{2.11}
\end{equation*}
\begin{equation*}
\int_1^\infty {x^{ - \textstyle{{\varepsilon + 1} \over p} - 1}[\int_0^{1 /
x} {k_\lambda (1,u)u^{r - \textstyle{\varepsilon \over p} - 1}du} ]dx} = O_3
(1),
\tag{2.12}
\label{2.12}
\end{equation*}
then for $\varepsilon \to 0^ + $，we have
\begin{equation*}
I_1 = \frac{1}{\varepsilon }(k_\lambda (r) + o_1 (1)) - O_1 (1),
\tag{2.13}
\label{2.13}
\end{equation*}
\begin{equation*}
I_2 = \frac{q}{1 + \varepsilon }(k_\lambda (r) + o_2 (1)) - O_2 (1),
\tag{2.14}
\label{2.14}
\end{equation*}
\begin{equation*}
I_3 = \frac{p}{1 + \varepsilon }(k_\lambda (r) + o_3 (1)) - O_3 (1).
\tag{2.15}
\label{2.15}
\end{equation*}
\end{lem}

{\bf Proof.} setting $u = y / x$,we have
\begin{equation*}
\begin{split}
 &I_1 = \int_1^\infty {\int_1^\infty {k_\lambda (x,y)x^{r -
\textstyle{\varepsilon \over p} - 1}y^{s - \textstyle{\varepsilon \over q} -
1}dxdy} } = \int_1^\infty {x^{ - 1 - \varepsilon }[\int_{1 / x}^\infty
{k_\lambda (1,u)u^{s - \textstyle{\varepsilon \over q} - 1}du} ]dx} \\
 &= \int_1^\infty {x^{ - 1 - \varepsilon }[\int_0^\infty {k_\lambda (1,u)u^{s
- \textstyle{\varepsilon \over q} - 1}du} ]dx} - \int_1^\infty {x^{ - 1 -
\varepsilon }[\int_0^{1 / x} {k_\lambda (1,u)u^{s - \textstyle{\varepsilon
\over q} - 1}du} ]dx} \\
 &= \frac{1}{\varepsilon }(k_\lambda (r) + o_1 (1)) - O_1 (1). \\
 \end{split}
\notag
\end{equation*}
Similarity we can prove (2.14) and (2.15), The lemma is proved.
\section{main results}

\begin{thm}\label{Theorem 3.1.} Let $p > 1$,$1 / p + 1 / q = 1$, $r + s = \lambda $, $f(x)$,
$g(x) \ge 0$, $\lambda > 0$

$F(x) = \int_0^x {f(t)dt} ,
G(x) = \int_0^x {g(t)dt} ,$ if $0 < \int_0^\infty {f^p(x)dx} < \infty $，

$0 < \int_0^\infty {g^q(x)dx} < \infty ,
k_\lambda (r): = \int\limits_0^\infty {k(u,1)u^{r - 1}du} $ is a positive
number, then the following two integral inequalities holds:
\begin{equation*}
\int_0^\infty {\int_0^\infty {x^{r - \textstyle{1 \over q} - 1}y^{s -
\textstyle{1 \over p} - 1}k_\lambda (x,y)F(x)G(y)dxdy} } < pqk_\lambda
(r)[\int_0^\infty {f^p(x)dx} ]^{\textstyle{1 \over p}}[\int_0^\infty
{g^q(x)dx} ]^{\textstyle{1 \over q}},
\tag{3.1}
\label{3.1}
\end{equation*}
And
\begin{equation*}
\int_0^\infty {[\int_0^\infty {x^{r - \textstyle{1 \over q} - 1}y^{s -
\textstyle{1 \over p}}k_\lambda (x,y)F(x)dx]^pdy} } < [qk_\lambda
(r)]^p\int_0^\infty {f^p(x)} dx,
\tag{3.2}
\label{3.2}
\end{equation*}
where the constant factors $pqk_\lambda (r)$ and $[qk_\lambda (r)]^p$ are
the best possible.\\
\end{thm}
{\bf Proof.}By H\"{o}lder's inequality with weight (cf. Kuang [9]) and Lemma
2.5, we have
\begin{equation*}
\begin{split}
 &\int_0^\infty {\int_0^\infty {x^{r - \textstyle{1 \over q} - 1}y^{s -
\textstyle{1 \over p} - 1}k_\lambda (x,y)F(x)G(y)} dxdy} \\
 &= \int_0^\infty {\int_0^\infty {k_\lambda (x,y)(y^{\textstyle{{s - 1} \over
p}}x^{\textstyle{r \over p} - 1}F(x))(x^{\textstyle{{r - 1} \over
q}}y^{\textstyle{s \over q} - 1}G(y))} dxdy} \\
 &\le \{\int_0^\infty {\int_0^\infty {k_\lambda (x,y)(y^{s - 1}x^{r -
p}F^p(x))} dxdy} \}^{\textstyle{1 \over p}} \\
 &\times \int_0^\infty {\int_0^\infty {k_\lambda (x,y)x^{r - 1}y^{s -
q}G^q(y))} dxdy} \\
 &= \{\int_0^\infty {\omega _\lambda (s,x)\left( {\frac{F(x)}{x}}
\right)^pdx} \}^{\frac{1}{p}}\{\int_0^\infty {\omega _\lambda (r,y)\left(
{\frac{G(y)}{y}} \right)^q} dy\}^{\frac{1}{q}} \\
 &= k_\lambda (r)\{\int_0^\infty {\left( {\frac{F(x)}{x}} \right)^pdx}
\}^{\frac{1}{p}}\{\int_0^\infty {\left( {\frac{G(y)}{y}} \right)^q}
dy\}^{\frac{1}{q}}, \\
 \end{split}
\notag
\end{equation*}
Then by Hardy inequality (3.1), (3.1) is valid.

By H\"{o}lder's inequality and Lemma 2.5, we get
\begin{equation*}
\begin{split}
& \int_0^\infty {x^{r - \textstyle{1 \over q} - 1}y^{s - \textstyle{1 \over
p}}k_\lambda (x,y)F(x)} dx \\
 &= \int_0^\infty {k_\lambda (x,y)(y^{\textstyle{{s - 1} \over
p}}x^{\textstyle{r \over p} - 1}F(x))(x^{\textstyle{{r - 1} \over
q}}y^{\textstyle{s \over q}})} dx \\
 &\le \{\int_0^\infty {k_\lambda (x,y)(y^{s - 1}x^{r - p}F^p(x))}
dx\}^{\textstyle{1 \over p}} \\
 &\times \{\int_0^\infty {k_\lambda (x,y)x^{r - 1}y^s} dx\}^{\textstyle{1
\over q}} \\
 &= k_\lambda ^{\textstyle{1 \over q}} (r)\{\int_0^\infty {k_\lambda
(x,y)(y^{s - 1}x^{r - p}F^p(x))} dx\}^{\textstyle{1 \over p}}, \\
 \end{split}
\notag
\end{equation*}
Hence, again applying Lemma 2.5, we obtain

\begin{equation*}
\begin{split}
 &\int_0^\infty {[\int_0^\infty {x^{r - \textstyle{1 \over q} - 1}y^{s -
\textstyle{1 \over p}}k_\lambda (x,y)F(x)} dx]^pdy} \\
 &\le k_\lambda ^{\textstyle{p \over q}} (r)\int_0^\infty {k_\lambda
(x,y)(y^{s - 1}x^{r - p}F^p(x))} dx \\
 &= k_\lambda ^p (r)\int_0^\infty {(\frac{F(x)}{x})^p} dx, \\
 \end{split}
\notag
\end{equation*}
then by Hardy inequality (2.1), (3.2) is valid.
\begin{equation*}
\begin{split}
 &\int_0^\infty {\int_0^\infty {x^{r - \textstyle{1 \over q} - 1}y^{s -
\textstyle{1 \over p} - 1}k_\lambda (x,y)F(x)G(y)} dxdy} \\
 &= \int_0^\infty {\int_0^\infty {x^{r - \textstyle{1 \over q} - 1}y^{s -
\textstyle{1 \over p}}k_\lambda (x,y)F(x)} dx(\frac{G(y)}{y})dy} \\
 &\le \{\int_0^\infty {[\int_0^\infty {x^{r - \textstyle{1 \over q} - 1}y^{s
- \textstyle{1 \over p}}k_\lambda (x,y)F(x)} dx]^pdy} \}^{\textstyle{1 \over
p}}\{\int_0^\infty {(\frac{G(y)}{y})^qdy} \}^{\textstyle{1 \over q}}, \\
 \end{split}
\notag
\end{equation*}
By (2.1) and (3.2), we have (3.1), Hence (3.2) and (3.1) are equivalent.
If the constant factor $pqk_\lambda (r)$ is not the best possible, then
there exists a positive constant $K$with$K < pqk_\lambda (r)$, thus (\ref{eq16}) is
still valid if we replace $pqk_\lambda (r)$ by $K$.

For sufficiently small $\varepsilon > 0$,Setting $f_\varepsilon (x)$，
$g_\varepsilon (x)$，$F_\varepsilon (x)$and$G_\varepsilon (x)$ as follow
\[
f_\varepsilon (x) = \left\{ {{\begin{array}{*{20}c}
 0 \hfill & {x \in (0,1)} \hfill \\
 {x^{ - \textstyle{1 \over p} - \textstyle{\varepsilon \over p}}} \hfill &
{x \in [1,\infty )} \hfill \\
\end{array} }} \right. \quad ,
\quad
g_\varepsilon (y) = \left\{ {{\begin{array}{*{20}c}
 0 \hfill & {y \in (0,1)} \hfill \\
 {y^{ - \textstyle{1 \over q} - \textstyle{\varepsilon \over q}}} \hfill &
{y \in [1,\infty )} \hfill \\
\end{array} }} \right.
\]
\[
F_\varepsilon (x) = \left\{ {{\begin{array}{*{20}c}
 {0,} \hfill & {x \in (0,1)} \hfill \\
 {\textstyle{q \over {1 - \varepsilon (q - 1)}}(x^{\textstyle{1 \over q} -
\textstyle{\varepsilon \over p}} - 1),} \hfill & {x \in [1,\infty )} \hfill
\\
\end{array} }} \right.,
\]

\[
G_\varepsilon (y) = \left\{ {{\begin{array}{*{20}c}
 {0,} \hfill & {y \in (0,1)} \hfill \\
 {\textstyle{p \over {1 - \varepsilon (p - 1)}}(y^{\textstyle{1 \over p} -
\textstyle{\varepsilon \over q}} - 1),} \hfill & {y \in [1,\infty )} \hfill
\\
\end{array} }} \right.,
\]
Let $\varphi (\varepsilon ) = \textstyle{{pq} \over {(1 - \varepsilon (q -
1))(1 - \varepsilon (p - 1))}}$, then $\varphi (\varepsilon ) \to pq$, as
$\varepsilon \to 0^ + $ and
\begin{equation*}
\{\int_0^\infty {f_\varepsilon ^p (x)dx} \}^{\textstyle{1 \over
p}}\{\int_0^\infty {g_\varepsilon ^q (x)dx} \}^{\textstyle{1 \over q}} =
\frac{1}{\varepsilon },
\tag{3.3}
\label{3.3}
\end{equation*}
\[
F_\varepsilon (x)G_\varepsilon (y) > \varphi (\varepsilon )(x^{\textstyle{1
\over q} - \textstyle{\varepsilon \over p}}y^{\textstyle{1 \over p} -
\textstyle{\varepsilon \over q}} - y^{\textstyle{1 \over p} -
\textstyle{\varepsilon \over q}} - x^{\textstyle{1 \over q} -
\textstyle{\varepsilon \over p}}),
\]
Hence
\begin{equation*}
\begin{split}
 &\int_0^\infty {\int_0^\infty {x^{r - \textstyle{1 \over q} - 1}y^{s -
\textstyle{1 \over p} - 1}k_\lambda (x,y)F_\varepsilon (x)G_\varepsilon
(y)dxdy} } \\
 &> \varphi (\varepsilon )\int_0^\infty {\int_0^\infty {k_\lambda (x,y)[x^{r -
\textstyle{\varepsilon \over p} - 1}y^{s - \textstyle{\varepsilon \over q} -
1} -x^{r - \textstyle{1 \over q} - 1}y^{s -
\textstyle{\varepsilon \over q} - 1} - x^{r - \textstyle{\varepsilon \over
p} - 1}y^{s - \textstyle{1 \over p} - 1}]dxdy} }, \\
 \end{split}
 \notag
\end{equation*}
By Lemma2.5,we obtian
\begin{equation*}
\begin{split}
 &\int_0^\infty {\int_0^\infty {x^{r - \textstyle{1 \over q} - 1}y^{s -
\textstyle{1 \over p} - 1}k_\lambda (x,y)F_\varepsilon (x)G_\varepsilon
(y)dxdy} } \\
 &> \varphi (\varepsilon )[\int_0^\infty {\int_0^\infty {k_\lambda (x,y)(x^{r -
\textstyle{\varepsilon \over p} - 1}y^{s - \textstyle{\varepsilon \over q} -
1}-x^{r - \textstyle{1 \over q} - 1}y^{s -
\textstyle{\varepsilon \over q} - 1})dxdy} } \\
&-\int_0^\infty {\int_0^\infty {k_\lambda(x,y) x^{r - \textstyle{\varepsilon \over
p} - 1}y^{s - \textstyle{1 \over p} - 1}dxdy} }] \\
 &= \varphi (\varepsilon )[\frac{1}{\varepsilon }(k_\lambda (r) + o_1 (1)) -
\frac{q}{1 + \varepsilon }(k_\lambda (r) + o_2 (1)) - \frac{p}{1 +
\varepsilon }(k_\lambda (r) + o_3 (1)) - O(1)], \\
 \end{split}
 \tag{3.4}\label{3.4}
\end{equation*}
If the constant factor$pqk_\lambda (r)$ in (3.1) is not the best possible,
then there exists a positive constant $K$, such that $K < pqk_\lambda (r)$
and (3.1) still remains valid if $pqk_\lambda (r)$ is replaced by $K$. In
particular by (3.2) and (3.3), we have
\begin{equation*}
\begin{split}
 &\varphi (\varepsilon )[k_\lambda (r) + o_1 (1) - \frac{\varepsilon q}{1 +
\varepsilon }(k_\lambda (r) + o_2 (1)) - \frac{\varepsilon p}{1 +
\varepsilon }(k_\lambda (r) + o_3 (1)) - \varepsilon O(1)] \\
 &< \varepsilon \int_0^\infty {\int_0^\infty {x^{r - \textstyle{1 \over q} -
1}y^{s - \textstyle{1 \over p} - 1}k_\lambda (x,y)F_\varepsilon
(x)G_\varepsilon (y)dxdy} } \\
 &< \varepsilon K\{\int_0^\infty {f_\varepsilon ^p (x)dx} \}^{\textstyle{1
\over p}}\{\int_0^\infty {g_\varepsilon ^q (x)dx} \}^{\textstyle{1 \over q}}
= K, \\
 \end{split}
\notag
\end{equation*}
Then $pqk_\lambda (r) \le K$ as $\varepsilon \to 0^ + $.This contradiction
shows that the constant factor $pqk_\lambda (r)$ in (3.1) is the best
possible.

If the constant factor $[qk_\lambda (r)]^p$ in (3.2) is not the best
possible, then there exists a positive constant $\tilde {K}$such that
$\tilde {K} < [qk_\lambda (r)]^p$and (3.2) still remains valid if
$[qk_\lambda (r)]^p$ is replaced by $\tilde {K}^p$. Then by Holder
inequality, (3.2) and Hardy inequality (2.1), we obtain

\begin{equation*}
\begin{split}
 &\int_0^\infty {\int_0^\infty {x^{r - \textstyle{1 \over q} - 1}y^{s -
\textstyle{1 \over p} - 1}k_\lambda (x,y)F(x)G(y)dxdy} } \\
 &= \int_0^\infty {\int_0^\infty {x^{r - \textstyle{1 \over q} - 1}y^{s -
\textstyle{1 \over p}}k_\lambda (x,y)F(x)\frac{G(y)}{y}dxdy} } \\
 &\le \{\int_0^\infty {(\int_0^\infty {x^{r - \textstyle{1 \over q} - 1}y^{s
- \textstyle{1 \over p}}k_\lambda (x,y)F(x)dx} )^pdy} \}^{\textstyle{1 \over
p}}\{\int_0^\infty {(\frac{G(y)}{y})^qdy} \}^{\textstyle{1 \over q}} \\
 &< p\tilde {K}[\int_0^\infty {f^p(x)dx} ]^{\textstyle{1 \over
p}}[\int_0^\infty {g^q(y)dy} ]^{\textstyle{1 \over q}}, \\
 \end{split}
 \notag
\end{equation*}
which gives that the constant factor $pqk_\lambda (r)$ in (3.1) is not the
best possible. This contradiction shows that the constant factor
$[qk_\lambda (r)]^p$ in (3.2) is the best possible. This proves the theorem.

\begin{thm}\label{Theorem 3.2.} Let $p > 1$,$1 / p + 1 / q = 1$,$p = \lambda - \alpha - 1 >
1$,$p = \lambda - \beta - 1 > 1$,$\alpha ,\beta > - 1,
f(x),
\quad
g(x) \ge 0,
F(x) = \int_0^x {f(t)dt} ,
G(x) = \int_0^x {g(t)dt} ,$ if $0 < \int_0^\infty {f^p(x)dx} < \infty $, $0 <
\int_0^\infty {g^q(x)dx} < \infty $, $k_\lambda (\alpha ): =
\int\limits_0^\infty {k(1,u)u^\alpha du} $ and $k_\lambda (\beta ): =
\int\limits_0^\infty {k(u,1)u^\beta du} $ are positive number ，then the
following two integral inequalities holds:
\begin{equation*}
\int_0^\infty {\int_0^\infty {k_\lambda (x,y)x^{\textstyle{\beta \over
q}}y^{\textstyle{\alpha \over p}}F(x)G(y)dxdy} } < pqk_\lambda ^{1 / p}
(\alpha )k_\lambda ^{1 / q} (\beta )[\int_0^\infty {f^p(x)dx} ]^{1 /
p}[\int_0^\infty {g^q(x)dx} ]^{\frac{1}{q}},
\tag{3.5}\label{3.5}
\end{equation*}
\begin{equation*}
\int_0^\infty {[\int_0^\infty {k_\lambda (x,y)x^{\textstyle{\beta \over
q}}y^{\textstyle{\alpha \over p} + 1}F(x)} dx]^pdy} < q^p[k_\lambda (\alpha
)]^{p - 1}k_\lambda (\beta )[\int_0^\infty {f^p(x)dx} ]^{1 /
p}[\int_0^\infty {g^q(x)dx} ]^{\frac{1}{q}},
\tag{3.6}
\label{3.6}
\end{equation*}
\end{thm}
The proof of Theorem 3.2 is similar to that of Theorem 3.1, so we omit it.

\begin{thm}\label{Theorem 3.3.} Let $p > 1$,$1 / p + 1 / q = 1$, $f(x)$, $g(x) \ge 0$, $\lambda> 0$, $F(x) = \int_x^\infty {f(t)dt} ,
G(x) = \int_x^\infty {g(t)dt} ,$if $0 < \int_0^\infty {(xf(x))^pdx} < \infty
$, $0 < \int_0^\infty {(xg(x))^qdx} < \infty ,
k_\lambda (p): = \int\limits_0^\infty {k(u,1)u^{\textstyle{1 \over p} -
1}du} $ is a positive number, then the following two integral inequalities
holds:
\begin{equation*}
\int_0^\infty {\int_0^\infty {k_\lambda (x,y)F(x)G(y)dxdy} } < pqk_\lambda
(r)[\int_0^\infty {(xf(x))^pdx} ]^{\textstyle{1 \over p}}[\int_0^\infty
{(xg(x))^qdx} ]^{\textstyle{1 \over q}},
\tag{3.7}
\label{3.7}
\end{equation*}
And
\begin{equation*}
\int_0^\infty {[\int_0^\infty {k_\lambda (x,y)F(x)dx]^pdy} } < [qk_\lambda
(r)]^p\int_0^\infty {(xf(x))^pdx},
\tag{3.8}
\label{3.8}
\end{equation*}
where the constant factors $pqk_\lambda (r)$ and $[qk_\lambda (r)]^p$ are
the best possible.
\end{thm}
The proof of Theorem 3.3 is similar to that of Theorem 3.1, so we omit it.
\begin{thm}\label{Theorem 3.4.} Let $0 < p < 1$,$1 / p + 1 / q = 1$, $r + s = \lambda
$,$f(x)$, $g(x) \ge 0$, $\lambda > 0$, $F(x) = \int_x^\infty {f(t)dt} ,
G(x) = \int_x^\infty {g(t)dt} , $if $0 < \int_0^\infty {f^p(x)dx} < \infty $,

$0 < \int_0^\infty {g^q(x)dx} < \infty ,
k_\lambda (r): = \int\limits_0^\infty {k(u,1)u^{r - 1}du} $ is a positive
number, then the following two integral inequalities holds:
\begin{equation*}
\int_0^\infty {\int_0^\infty {x^{r - \textstyle{1 \over q} - 1}y^{s -
\textstyle{1 \over p} - 1}k_\lambda (x,y)F(x)G(y)dxdy} } > ( - pqk_\lambda
(r))[\int_0^\infty {f^p(x)dx} ]^{\textstyle{1 \over p}}[\int_0^\infty
{g^q(x)dx} ]^{\textstyle{1 \over q}},
\tag{3.9}
\label{3.9}
\end{equation*}
And
\begin{equation*}
\int_0^\infty {[\int_0^\infty {x^{r - \textstyle{1 \over q} - 1}y^{s -
\textstyle{1 \over p}}k_\lambda (x,y)F(x)dx]^pdy} } > [ - qk_\lambda
(r)]^p\int_0^\infty {f^p(x)} dx,
\tag{3.10}
\label{3.10}
\end{equation*}
where the constant factors $[ - pqk_\lambda (r)]$ and $[ - qk_\lambda
(r)]^p$ are the best possible.
\end{thm}
\begin{thm}\label{Theorem 3.5.} Let $0 < p < 1$,$1 / p + 1 / q = 1$, $f(x)$, $g(x) \ge 0$,
$\lambda > 0$, $F(x) = \int_x^\infty {f(t)dt} ,
G(x) = \int_x^\infty {g(t)dt} , $if $0 < \int_0^\infty {(xf(x))^pdx} < \infty
$, $0 < \int_0^\infty {(xg(x))^qdx} < \infty ,
k_\lambda (p): = \int\limits_0^\infty {k(u,1)u^{\textstyle{1 \over p} -
1}du} $ is a positive number, then the following two integral inequalities
holds:
\begin{equation*}
\int_0^\infty {\int_0^\infty {k_\lambda (x,y)F(x)G(y)dxdy} } > [ -
pqk_\lambda (p)][\int_0^\infty {(xf(x))^pdx} ]^{\textstyle{1 \over
p}}[\int_0^\infty {(xg(x))^qdx} ]^{\textstyle{1 \over q}},
\tag{3.11}
\label{3.11}
\end{equation*}
And
\begin{equation*}
\int_0^\infty {[\int_0^\infty {k_\lambda (x,y)F(x)dx]^pdy} } > [ -
qk_\lambda (p)]^p\int_0^\infty {(xf(x))^pdx},
\tag{3.12}
\label{3.12}
\end{equation*}
where the constant factors $[ - pqk_\lambda (r)]$ and $[ - qk_\lambda
(r)]^p$ are the best possible.
\end{thm}
\section{Discrete analogous}
\begin{thm}\label{Theorem 4.1.} Let $p > 1$,$1 / p + 1 / q = 1$, $r + s = \lambda $,$a_n
$, $b_n \geq 0$, $A_n = \sum\limits_{k = 1}^n {a_k } $, $B_n = \sum\limits_{k =
1}^n {b_k } $, if $k(u,1)u^{r - 1}$ and $k(1,u)u^{s - 1}$ are decreasing in
$(0,\infty )$ and strictly decreasing in a subinterval of $(0,\infty )$, $0
< \sum\limits_{n = 1}^\infty {a_n^p } < \infty $, $0 < \sum\limits_{n =
1}^\infty {b_n^q } < \infty $, then the following two inequalities holds:
\begin{equation*}
\sum\limits_{n = 1}^\infty {\sum\limits_{m = 1}^\infty {m^{r - \textstyle{1
\over q} - 1}n^{s - \textstyle{1 \over p} - 1}k_\lambda (m,n)A_m B_n } } <
pqk_\lambda (r)[\sum\limits_{n = 1}^\infty {a_n^p } ]^{\textstyle{1 \over
p}}[\sum\limits_{n = 1}^\infty {b_n^q } ]^{\textstyle{1 \over q}},
\tag{4.1}
\label{4.1}
\end{equation*}
And
\begin{equation*}
\sum\limits_{n = 1}^\infty {[\sum\limits_{m = 1}^\infty {m^{r - \textstyle{1
\over q} - 1}n^{s - \textstyle{1 \over p}}k_\lambda (m,n)A_m } } ]^p <
[qk_\lambda (r)]^p\sum\limits_{n = 1}^\infty {a_n^p },
\tag{4.2}
\label{4.2}
\end{equation*}
where the constant factors $pqk_\lambda (r)$ and $[qk_\lambda (r)]^p$ are
the best possible.
\end{thm}
The proof of Theorem 4.1 is similar to that of Theorem 3.1, so we omit it.
\begin{thm}\label{
Theorem 4.2.} Let $p > 1$,$1 / p + 1 / q = 1$, $p = \lambda - \alpha - 1 >
1$,$p = \lambda - \beta - 1 > 1$, $\alpha ,\beta > - 1,
a_n ,
b_n \geq 0 ,
A_n = \sum\limits_{k = 1}^n {a_k } ,
B_n = \sum\limits_{k = 1}^n {b_k } , $if $k(u,1)u^\alpha $ and $k(1,u)u^\beta
$ are decreasing in $(0,\infty )$ and strictly decreasing in a subinterval of
$(0,\infty )$, $0 < \sum\limits_{n = 1}^\infty {a_n^p } < \infty $, $0 <
\sum\limits_{n = 1}^\infty {b_n^q } < \infty $, then the following two
inequalities holds:
\begin{equation*}
\sum\limits_{n = 1}^\infty {\sum\limits_{m = 1}^\infty {m^{\textstyle{\beta
\over q}}n^{\textstyle{\alpha \over p}}k_\lambda (m,n)A_m B_n } } <
pqk_\lambda (\alpha )[\sum\limits_{n = 1}^\infty {a_n^p } ]^{\textstyle{1
\over p}}[\sum\limits_{n = 1}^\infty {b_n^q } ]^{\textstyle{1 \over q}},
\tag{4.3}
\label{4.3}
\end{equation*}
And
\begin{equation*}
\sum\limits_{n = 1}^\infty {[\sum\limits_{m = 1}^\infty {m^{\textstyle{\beta
\over q}}n^{\textstyle{\alpha \over p} + 1}k_\lambda (m,n)A_m } } ]^p <
q^p[k_\lambda (\alpha )]^{p - 1}k_\lambda (\beta )\sum\limits_{n = 1}^\infty
{a_n^p },
\tag{4.4}
\label{4.4}
\end{equation*}
\end{thm}

\begin{thm}\label{Theorem 4.3.} Let $p > 1$,$1 / p + 1 / q = 1$, $\lambda > 0$, $a_n $, $b_n \geq 0
$, $A_n = \sum\limits_{k = 1}^n {a_k } $, $B_n = \sum\limits_{k = 1}^n {b_k
} $, $k_\lambda (p): = \int\limits_0^\infty {k(u,1)u^{\textstyle{1 \over p}
- 1}du} $,if $k(u,1)u^{\textstyle{1 \over p} - 1}$ and
$k(1,u)u^{\textstyle{1 \over q} - 1}$ are decreasing in $(0,\infty )$ and
strictly decreasing in a subinterval of $(0,\infty )$, $0 < \sum\limits_{n =
1}^\infty {a_n^p } < \infty $, $0 < \sum\limits_{n = 1}^\infty {b_n^q } <
\infty $, then the following two inequalities holds:
\begin{equation*}
\sum\limits_{n = 1}^\infty {\sum\limits_{m = 1}^\infty {k_\lambda (m,n)A_m
B_n } } < pqk_\lambda (p)[\sum\limits_{n = 1}^\infty {(na_n )^p}
]^{\textstyle{1 \over p}}[\sum\limits_{n = 1}^\infty {(nb_n )^p}
]^{\textstyle{1 \over q}},
\tag{4.5}
\label{4.5}
\end{equation*}
And
\begin{equation*}
\sum\limits_{n = 1}^\infty {[\sum\limits_{m = 1}^\infty {k_\lambda (m,n)A_m
} } ]^p < [qk_\lambda (p)]^p\sum\limits_{n = 1}^\infty {(na_n )^p},
\tag{4.6}
\label{4.6}
\end{equation*}
where the constant factors $pqk_\lambda (p)$ and $[qk_\lambda (p)]^p$ are
the best possible.
\end{thm}
\begin{thm}\label{Theorem 4.4.} Let $0 < p < 1$,$1 / p + 1 / q = 1$, $r + s = \lambda $, $a_n
$, $b_n \geq 0 $, $A_n = \sum\limits_{k = 1}^n {a_k } $, $B_n = \sum\limits_{k =
1}^n {b_k } $, if $k(u,1)u^{r - 1}$ and $k(1,u)u^{s - 1}$ are decreasing in
$(0,\infty )$ and strictly decreasing in a subinterval of $(0,\infty )$, $0
< \sum\limits_{n = 1}^\infty {a_n^p } < \infty $, $0 < \sum\limits_{n =
1}^\infty {b_n^q } < \infty $, then the following two inequalities holds:
\begin{equation*}
\sum\limits_{n = 1}^\infty {\sum\limits_{m = 1}^\infty {m^{r - \textstyle{1
\over q} - 1}n^{s - \textstyle{1 \over p} - 1}k_\lambda (m,n)A_m B_n } } > [
- pqk_\lambda (r)][\sum\limits_{n = 1}^\infty {a_n^p } ]^{\textstyle{1 \over
p}}[\sum\limits_{n = 1}^\infty {b_n^q } ]^{\textstyle{1 \over q}},
\tag{4.7}
\label{4.7}
\end{equation*}
And
\begin{equation*}
\sum\limits_{n = 1}^\infty {[\sum\limits_{m = 1}^\infty {m^{r - \textstyle{1
\over q} - 1}n^{s - \textstyle{1 \over p}}k_\lambda (m,n)A_m } } ]^p > [ -
qk_\lambda (r)]^p\sum\limits_{n = 1}^\infty {a_n^p },
\tag{4.8}
\label{4.8}
\end{equation*}
where the constant factors $[ - pqk_\lambda (r)]$ and $[ - qk_\lambda
(r)]^p$ are the best possible.
\end{thm}

\begin{thm}\label{Theorem 4.5.} Let $0 < p < 1$,$1 / p + 1 / q = 1$, $\lambda > 0$, $a_n
$, $b_n \geq 0 $, $A_n = \sum\limits_{k = 1}^n {a_k } $, $B_n = \sum\limits_{k =
1}^n {b_k } $, $k_\lambda (p): = \int\limits_0^\infty {k(u,1)u^{\textstyle{1
\over p} - 1}du} $, if $k(u,1)u^{\textstyle{1 \over p} - 1}$ and
$k(1,u)u^{\textstyle{1 \over q} - 1}$ are decreasing in $(0,\infty )$ and
strictly decreasing in a subinterval of $(0,\infty )$, $0 < \sum\limits_{n =
1}^\infty {a_n^p } < \infty $, $0 < \sum\limits_{n = 1}^\infty {b_n^q } <
\infty $, then the following two inequalities holds:
\begin{equation*}
\sum\limits_{n = 1}^\infty {\sum\limits_{m = 1}^\infty {k_\lambda (m,n)A_m
B_n } } > [ - pqk_\lambda (p)][\sum\limits_{n = 1}^\infty {(na_n )^p}
]^{\textstyle{1 \over p}}[\sum\limits_{n = 1}^\infty {(nb_n )^p}
]^{\textstyle{1 \over q}},
\tag{4.9}
\label{4.9}
\end{equation*}
And
\begin{equation*}
\sum\limits_{n = 1}^\infty {[\sum\limits_{m = 1}^\infty {k_\lambda (m,n)A_m
} } ]^p > [ - qk_\lambda (p)]^p\sum\limits_{n = 1}^\infty {(na_n )^p},
\tag{4.10}
\label{4.10}
\end{equation*}
where the constant factors $[ - pqk_\lambda (p)]$ and $[ - qk_\lambda
(p)]^p$ are the best possible.
\end{thm}

\section{some particular results}

(1) $k_\lambda (x,y) = \frac{1}{\vert x - y\vert ^\lambda }$, by
Lemma 2.3， we have
\[
k_\lambda (r): = \int\limits_0^\infty {\frac{u^{r - 1}}{\vert 1 - u\vert
^\lambda }du} = B(r,1 - \lambda ) + B(s,1 - \lambda ),
\]
By Theorem 3.1 and 4.1, we have
\begin{equation*}
\int_0^\infty {\int_0^\infty {\frac{x^{r - \textstyle{1 \over q} - 1}y^{s -
\textstyle{1 \over p} - 1}F(x)G(y)}{\vert x - y\vert ^\lambda }dxdy} } <
pq[B(r,1 - \lambda ) + B(s,1 - \lambda )][\int_0^\infty {f^p(x)dx}
]^{\textstyle{1 \over p}}[\int_0^\infty {g^q(x)dx} ]^{\textstyle{1 \over
q}},
\tag{5.1}
\label{5.1}
\end{equation*}
\begin{equation*}
\int_0^\infty {[\int_0^\infty {\frac{x^{r - \textstyle{1 \over q} - 1}y^{s -
\textstyle{1 \over p}}F(x)}{\vert x - y\vert ^\lambda }dx]^pdy} } < [q(B(r,1
- \lambda ) + B(s,1 - \lambda ))]^p\int_0^\infty {f^p(x)} dx,
\tag{5.2}
\label{5.2}
\end{equation*}
\\ (2) $k_\lambda (x,y) = \frac{\ln (x / y)}{x^\lambda - y^\lambda }$ by Lemma
2.3， we have
\[
k_\lambda (r): = \int\limits_0^\infty {\frac{\ln uu^{r - 1}}{1 - u^\lambda
}du} = [\textstyle{\pi \over {\lambda \sin (r / \lambda )}}]^2
\]
By Theorem 3.1 and 4.1, we have
\begin{equation*}
\int_0^\infty {\int_0^\infty {\frac{\ln (x / y)x^{r - \textstyle{1 \over q}
- 1}y^{s - \textstyle{1 \over p} - 1}F(x)G(y)}{x^\lambda - y^\lambda }dxdy}
} < pq[\textstyle{\pi \over {\lambda \sin (r / \lambda )}}]^2[\int_0^\infty
{f^p(x)dx} ]^{\textstyle{1 \over p}}[\int_0^\infty {g^q(x)dx}
]^{\textstyle{1 \over q}},
\tag{5.3}
\label{5.3}
\end{equation*}
\begin{equation*}
\int_0^\infty {[\int_0^\infty {\frac{\ln (x / y)x^{r - \textstyle{1 \over q}
- 1}y^{s - \textstyle{1 \over p}}F(x)}{x^\lambda - y^\lambda }dx]^pdy} } <
[q[\textstyle{\pi \over {\lambda \sin (r / \lambda )}}]^2]^p\int_0^\infty
{f^p(x)} dx,
\tag{5.4}
\label{5.4}
\end{equation*}
\\ (3) $k_\lambda (x,y) = \frac{1}{\vert x - y\vert ^\beta (\max
\{x,y\})^{\lambda - \beta }}$, ($0 < \beta < 1)$ by Lemma 2.3， we have

\[
k_\lambda (r): = \int\limits_0^\infty {\frac{u^{r - 1}}{\vert 1 - u\vert
^\beta (\max \{1,u\})^{\lambda - \beta }}du} = B(r,1 - \beta ) + B(s,1 -
\beta ),
\]
By Theorem 3.1 and 4.1, we have
\begin{equation*}
\int_0^\infty {\int_0^\infty {\frac{x^{r - \textstyle{1 \over q} - 1}y^{s -
\textstyle{1 \over p} - 1}F(x)G(y)}{\vert x - y\vert ^\beta (\max
\{x,y\})^{\lambda - \beta }}dxdy} } < pq[B(r,1 - \beta ) + B(s,1 - \beta
)][\int_0^\infty {f^p(x)dx} ]^{\textstyle{1 \over p}}[\int_0^\infty
{g^q(x)dx} ]^{\textstyle{1 \over q}},
\tag{5.5}
\label{5.5}
\end{equation*}
\begin{equation*}
\int_0^\infty {[\int_0^\infty {\frac{x^{r - \textstyle{1 \over q} - 1}y^{s -
\textstyle{1 \over p}}F(x)}{\vert x - y\vert ^\beta (\max \{x,y\})^{\lambda
- \beta }}dx]^pdy} } < [q(B(r,1 - \beta ) + B(s,1 - \beta ))]^p\int_0^\infty
{f^p(x)} dx,
\tag{5.6}
\label{5.6}
\end{equation*}
\\ (4) $k_\lambda (x,y) = \frac{(\min \{x,y\})^{\beta - \lambda }}{\vert x -
y\vert ^\beta }$, ($0 < \beta < 1)$ by Lemma 2.3， we have
\[
k_\lambda (r): = \int\limits_0^\infty {\frac{(\min \{1,u\})^{\beta - \lambda
}u^{r - 1}}{\vert 1 - u\vert ^\beta }du} = B(\beta - r,1 - \beta ) + B(\beta
- s,1 - \beta ),
\]
By Theorem 3.1 and 4.1, we have
\begin{equation*}
\begin{split}
 &\int_0^\infty {\int_0^\infty {\frac{(\min \{x,y\})^{\beta - \lambda }x^{r -
\textstyle{1 \over q} - 1}y^{s - \textstyle{1 \over p} - 1}F(x)G(y)}{\vert x
- y\vert ^\beta }dxdy} } \\
 &< pq[(\beta - r,1 - \beta ) + B(\beta - s,1 - \beta )][\int_0^\infty
{f^p(x)dx} ]^{\textstyle{1 \over p}}[\int_0^\infty {g^q(x)dx}
]^{\textstyle{1 \over q}}, \\
\end{split}
\tag{5.7}
\label{5.7}
\end{equation*}
\begin{equation*}
\int_0^\infty {[\int_0^\infty {\frac{(\min \{x,y\})^{\beta - \lambda }x^{r -
\textstyle{1 \over q} - 1}y^{s - \textstyle{1 \over p}}F(x)}{\vert x -
y\vert ^\beta }dx]^pdy} } < [q((\beta - r,1 - \beta ) + B(\beta - s,1 -
\beta ))]^p\int_0^\infty {f^p(x)} dx,
\tag{5.8}
\label{5.8}
\end{equation*}
\\ (5) $k_\lambda (x,y) = \frac{\vert x^\beta - y^\beta \vert }{(\max
\{x,y\})^{\lambda + \beta }}$, ($\beta > - \min \{r,s\})$ by Lemma 2.3， we
have
\[
k_\lambda (r): = \int\limits_0^\infty {\frac{\vert 1 - u^\beta \vert u^{r -
1}}{(\max \{1,u\})^{\lambda + \beta }}du} = \frac{\vert \beta \vert (r(r +
\beta ) + s(s + \beta ))}{rs(r + \beta )(s + \beta )},
\]
By Theorem 3.1 and 4.1, we have
\begin{equation*}
\begin{split}
 &\int_0^\infty {\int_0^\infty {\frac{\vert x^\beta - y^\beta \vert x^{r -
\textstyle{1 \over q} - 1}y^{s - \textstyle{1 \over p} - 1}F(x)G(y)}{(\max
\{x,y\})^{\lambda + \beta }}dxdy} } \\
 &< pq\frac{\vert \beta \vert (r(r + \beta ) + s(s + \beta ))}{rs(r + \beta
)(s + \beta )}[\int_0^\infty {f^p(x)dx} ]^{\textstyle{1 \over
p}}[\int_0^\infty {g^q(x)dx} ]^{\textstyle{1 \over q}}, \\
 \end{split}
\tag{5.9}
\label{5.9}
\end{equation*}
\begin{equation*}
\int_0^\infty {[\int_0^\infty {\frac{\vert x^\beta - y^\beta \vert x^{r -
\textstyle{1 \over q} - 1}y^{s - \textstyle{1 \over p}}F(x)}{(\max
\{x,y\})^{\lambda + \beta }}dx]^pdy} } < [q\frac{\vert \beta \vert (r(r +
\beta ) + s(s + \beta ))}{rs(r + \beta )(s + \beta )}]^p\int_0^\infty
{f^p(x)} dx,
\tag{5.10}
\label{5.10}
\end{equation*}
\\ (6) $k_\lambda (x,y) = \frac{\vert \ln (x / y)\vert }{(\max \{x,y\})^\lambda
}$, ($0 < \beta < 1)$ by Lemma 2.3， we have
\[
k_\lambda (r): = \int\limits_0^\infty {\frac{\vert \ln u\vert u^{r -
1}}{(\max \{1,u\})^\lambda }du} = \frac{1}{r^2} + \frac{1}{s^2},
\]
By Theorem 3.1 and 4.1, we have
\begin{equation*}
\int_0^\infty {\int_0^\infty {\frac{\vert \ln (x / y)\vert x^{r -
\textstyle{1 \over q} - 1}y^{s - \textstyle{1 \over p} - 1}F(x)G(y)}{(\max
\{x,y\})^\lambda }dxdy} } < pq(\frac{1}{r^2} + \frac{1}{s^2})[\int_0^\infty
{f^p(x)dx} ]^{\textstyle{1 \over p}}[\int_0^\infty {g^q(x)dx}
]^{\textstyle{1 \over q}},
\tag{5.11}
\label{5.11}
\end{equation*}
\begin{equation*}
\int_0^\infty {[\int_0^\infty {\frac{\vert \ln (x / y)\vert x^{r -
\textstyle{1 \over q} - 1}y^{s - \textstyle{1 \over p}}F(x)}{(\max
\{x,y\})^\lambda }dx]^pdy} } < [q(\frac{1}{r^2} +
\frac{1}{s^2})]^p\int_0^\infty {f^p(x)} dx,
\tag{5.12}
\label{5.12}
\end{equation*}
\\ (7) $k_\lambda (x,y) = \frac{\vert \ln (x / y)\vert }{x^\lambda + y^\lambda
}$, ($0 < \beta < 1)$ by Lemma 2.3， we have
\[
k_\lambda (r): = \int\limits_0^\infty {\frac{\vert \ln u\vert u^{r - 1}}{1 +
u^\lambda }du} = \sum\limits_{n = 1}^\infty {( - 1)^n\textstyle{2 \over
{(\lambda n + r)^2}}},
\]
By Theorem 3.1 and 4.1, we have
\begin{equation*}
\begin{split}
 &\int_0^\infty {\int_0^\infty {\frac{\vert \ln (x / y)\vert x^{r -
\textstyle{1 \over q} - 1}y^{s - \textstyle{1 \over p} -
1}F(x)G(y)}{x^\lambda + y^\lambda }dxdy} } \\
 &< pq(\sum\limits_{n = 1}^\infty {( - 1)^n\textstyle{2 \over {(\lambda n +
r)^2}}} )^{\textstyle{1 \over p}}(\sum\limits_{n = 1}^\infty {( -
1)^n\textstyle{2 \over {(\lambda n + s)^2}}} )^{\textstyle{1 \over
q}}[\int_0^\infty {f^p(x)dx} ]^{\textstyle{1 \over p}}[\int_0^\infty
{g^q(x)dx} ]^{\textstyle{1 \over q}}, \\
\end{split}
\tag{5.13}
\label{5.13}
\end{equation*}
\begin{equation*}
\int_0^\infty {[\int_0^\infty {\frac{\vert \ln (x / y)\vert x^{r -
\textstyle{1 \over q} - 1}y^{s - \textstyle{1 \over p}}F(x)}{x^\lambda +
y^\lambda }dx]^pdy} } < q^p(\sum\limits_{n = 1}^\infty {( - 1)^n\textstyle{2
\over {(\lambda n + r)^2}}} )^{p - 1}\sum\limits_{n = 1}^\infty {( -
1)^n\textstyle{2 \over {(\lambda n + s)^2}}} \int_0^\infty {f^p(x)} dx,
\tag{5.14}
\label{5.14}
\end{equation*}

 \end{document}